\documentclass[11pt]{amsart}

\usepackage{amsmath}
\usepackage{eucal}
\usepackage{amscd}
\usepackage{wrapfig}
\usepackage{graphicx}
\usepackage{amssymb}
\usepackage{latexsym}
\usepackage{amsfonts}
\newtheorem{lemma}{Lemma}[section]
\newtheorem{thm}{Theorem}[section]
\newtheorem{prop}{Proposition}[section]
\newtheorem{cor}{Corollary}[section]
\newtheorem{rem}{Remark}[section]

\numberwithin{equation}{section}
\title{Generalized Riemann sums} 

\author{Toshikazu Sunada}

\address{
School of Interdisciplinary Mathematical Sciences, Meiji University, Nakano 4-21-1,  Nakano-ku, Tokyo, 164-8525 Japan}
\address{Advanced Institute for Materials Research, Tohoku University, 2-1-1 Katahira, Aoba-ku, Sendai, 980-8577 Japan}
\email{sunada@meiji.ac.jp}
\thanks{The author is partially supported by JSPS Grant-in-Aid for Scientific Research (No. 15H02055).}
\date{}
\keywords{constant density, coprime pairs, primitive Pythagorean triples, quasicrystal, rational points on the unit circle}

\begin{document}
\maketitle

\begin{abstract}
The primary aim of this chapter is, commemorating the 150th anniversary of Riemann's death, to explain how the idea of {\it Riemann sum} is linked to other branches of mathematics. The materials I treat are more or less classical and elementary, thus available to the ``common mathematician in the streets." However one may still see here interesting inter-connection and cohesiveness in mathematics.
\end{abstract}

\section{Introduction}

On Gauss's recommendation, Bernhard Riemann presented the paper {\it \"{U}ber die Darstellbarkeit einer Function durch eine trigonometrische Reihe} to the Council of G\"{o}ttingen University as his Habilitationsschrift at the first stage in December of 1853.\footnote{Habilitationsschrift is a thesis for qualification to become a lecturer. The famous lecture {\it \"{U}ber die Hypothesen welche der Geometrie zu Grunde liegen} delivered on 10 June 1854 was for the final stage of his Habilitationsschrift.} As the title clearly suggests, the aim of his essay was to lay the foundation for the theory of trigonometric series (Fourier series\index{Fourier series} in today's term).\footnote{The English translation is ``On the representability of a function by a trigonometric series". His essay was published only after his death in the Abhandlungen der K\"{o}niglichen Gesellschaft der Wissenschaften zu G\"{o}ttingen (Proceedings of the Royal Philosophical Society at G\"{o}ttingen), vol. 13, (1868), pages 87--132.} 

The record of previous work by other mathematicians, to which Riemann devoted three sections of the essay, tells us that the Fourier series had been used to {\it represent} general solutions of the wave equation and the heat equation without any convincing proof of convergence. For instance, Fourier claimed, in his study of the heat equation (1807, 1822), that if we put
\begin{equation}\label{eq:fourier1}
a_n=\frac{1}{\pi}\int_{-\pi}^{\pi}f(x)\sin nx~\!dx,\quad b_n=\frac{1}{\pi}\int_{-\pi}^{\pi}f(x)\cos nx~\!dx,
\end{equation}
then
\begin{equation}\label{eq:fourier}
f(x)=\frac{1}{2}b_0+(a_1\sin x+b_1\cos x)+(a_2\sin 2x+b_2\cos 2x)+\cdots
\end{equation}
without any restrictions on the function $f(x)$.   
But this is not true in general as is well known. What is worse (though, needless to say, the significance of his paper as a historical document cannot be denied) is his claim that the integral of an ``arbitrary" function is meaningful as the area under/above the associated graph.

L. Dirichlet, a predecessor of Riemann, was the first who gave a solid proof for convergence in a special case. Actually he proved that the right-hand side of (\ref{eq:fourier}) converges to $\displaystyle\frac{1}{2}\big(f(x+0)+f(x-0)\big)$ for a class of functions including piecewise monotone continuous functions (1829). Stimulated by Dirichlet's study, Riemann made considerable progress on the convergence problem. In the course of his discussion, he gave a precise notion of integrability of a function,\footnote{See Section 4 in his essay, entitled ``\"{U}ber der Begriff eines bestimmten Integrals und den Umfang seiner G\"{u}ltigkeit" (On the concept of a definite integral and the extent of its validity), pages 101-103.} and then obtained a condition for an integrable function to be representable by a Fourier series. Furthermore he proved that the Fourier coefficients for any integrable function $a_n, b_n$ converge to zero as $n\to\infty$. This theorem, which was generalized by Lebesgue later to a broader class of functions, is to be called the Riemann-Lebesgue theorem, and is of importance in Fourier analysis and asymptotic analysis.

\begin{figure}[htbp]
\includegraphics[width=.95\linewidth]{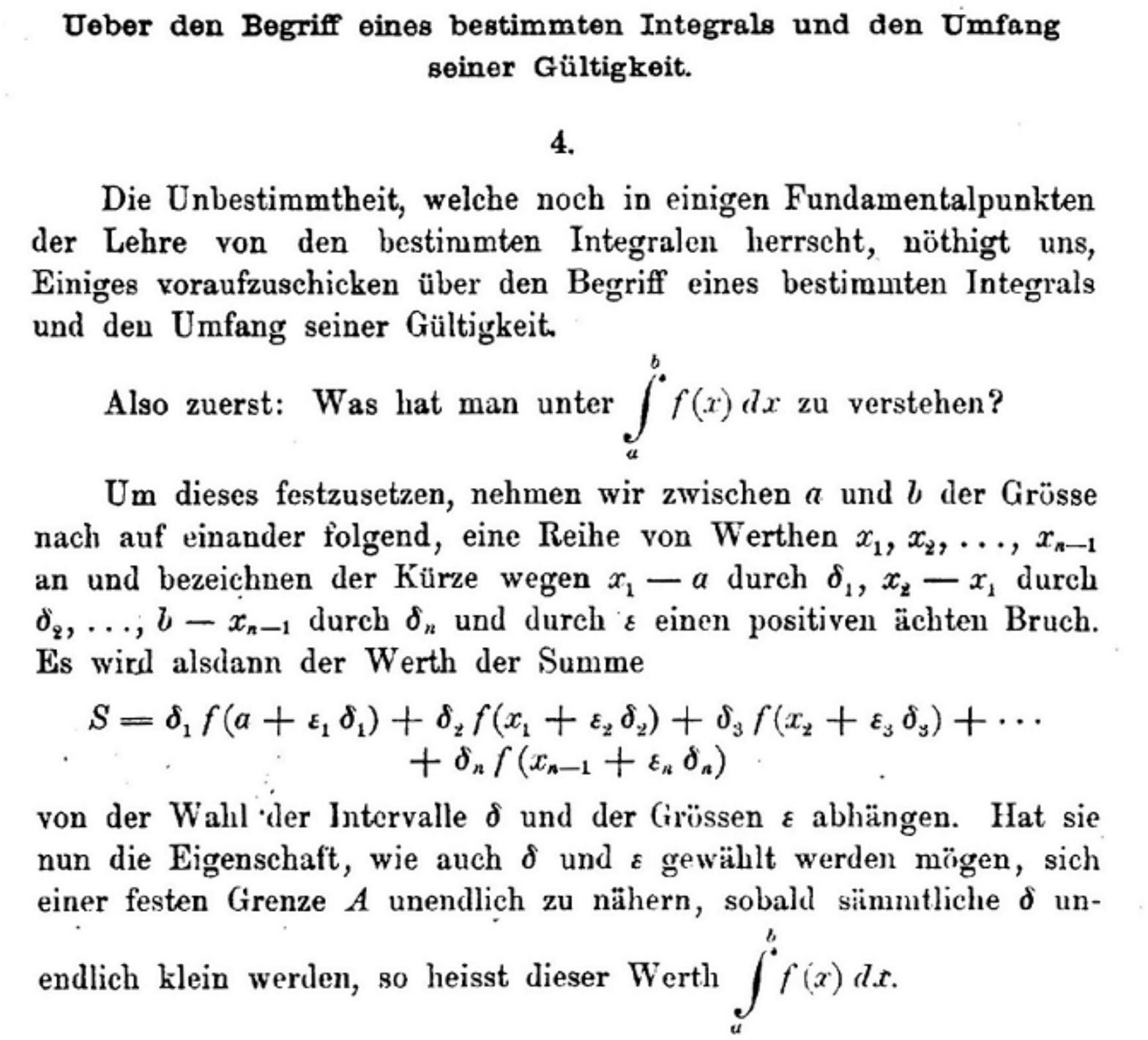}
\vspace{-7cm}
\caption{Riemann's paper}\label{fig:riemannpapereps}
\end{figure}

What plays a significant role in Riemann's definition of integrals is the notion of {\it Riemann sum}\index{Riemann sum}, which, if we use his notation (Fig.~\!\ref{fig:riemannpapereps}), is expressed as
$$
S=\delta_1f(a+\epsilon_1\delta_1)+\delta_2f(x_1+\epsilon_2\delta_2)+\delta_3f(x_3+\epsilon_3\delta_3)+\cdots+\delta_nf(x_{n-1}+\epsilon_n\delta_n).
$$

\noindent Here $f(x)$ is a function on the closed interval $[a,b]$,  $a=x_0<x_1<x_2<\cdots<x_{n-1}<x_n=b$, and $\delta_i=x_{i}-x_{i-1}$ ($i=1,2,\ldots,n$). If $S$ converges to $A$ as $\max_i \delta_i$ goes to $0$ whatever $\epsilon_i$ with $0<\epsilon_i<1$ ($i=1,\ldots,n$) are chosen (thus $x_{k-1}+\epsilon_k\delta_k\in [x_{k-1},x_k]$), then the value $A$ is written as $\displaystyle\int_a^bf(x)dx$, and $f(x)$ is called {\it Riemann integrable}\index{Riemann integrable}. For example, every continuous function is Riemann integrable as we learn in calculus.

Compared with Riemann's other supereminent works, his essay looks unglamorous. Indeed, from today's view, his formulation of integrability is no more than routine. But the harbinger must push forward through the total dark without any definite idea of the direction. All he needs is a torch of intelligence. 

The primary aim of this chapter is {\it not} to present the subsequent development after Riemann's work on integrals such as the contribution by C. Jordan (1892)\footnote{Jordan introduced a measure (Jordan measure) which fits in with Riemann integral. A bounded set is Jordan measurable if and only if its indicator function is Riemann integrable.}, G. Peano (1887), H. L. Lebesgue (1892), T. J. Stieltjes (1894), and K. Ito (1942)\footnote{Ito's integral (or stochastic integral) is a sort of generaization of Stieltjes integral. Stieltjes defined his integral $\displaystyle \int f(x)d\varphi(x)$ by means of a modified Riemann sum.}, but to explain how the idea of Riemann sum is linked to other branches of mathematics; for instance, some counting problems in elementary number theory and the theory of quasicrystals, the former having a long history and the latter being  an active field still in a state of flux.  

\medskip 

I am very grateful to Xueping Guang for drawing attention to Ref.~\!\cite{matt} which handles some notions closely related to the ones in the present chapter.

\vspace{-1cm}
\hspace{2cm}

\section{Generalized Riemann sums}

The notion of Riemann sum is immediately generalized to functions of several variables as follows. 

Let $\Delta=\{D_{\alpha}\}_{\alpha\in A}$ be a partition of $\mathbb{R}^d$
by  a countable family of 
bounded domains $D_{\alpha}$ with piecewise smooth boundaries satisfying

\smallskip

(i)~ 
${\rm mesh}(\Delta):=\displaystyle\sup_{\alpha\in A} d(D_{\alpha})<\infty$, where $d(D_{\alpha})$ is the diameter of $D_{\alpha}$,

\smallskip

(ii)~ there are only finitely many $\alpha$ such that $K\cap D_{\alpha}\neq \emptyset$ for any compact set $K\subset \mathbb{R}^d$.
\smallskip

\smallskip

We select a point ${\boldsymbol \xi}_{\alpha}$ from each $D_{\alpha}$, and put $\Gamma=\{\boldsymbol{\xi}_{\alpha}|~\alpha\in A\}$. The Riemann sum\index{Riemann sum} $\sigma(f,\Delta,\Gamma)$ for a function $f$ on $\mathbb{R}^d$ with compact support is defined by 
$$
\sigma(f,\Delta,\Gamma)=\sum_{\alpha\in A}f({\boldsymbol \xi}_{\alpha}){\rm vol}(D_{\alpha}),$$
where ${\rm vol}(D_{\alpha})$ is the volume of $D_{\alpha}$.  
Note that $f({\boldsymbol \xi}_{\alpha})=0$ for all but finitely many $\alpha$ because of Property (ii).

If the limit 
\begin{equation*}
\displaystyle\lim_{{\rm mesh}(\Delta)\to 0}\sigma(f,\Delta,\Gamma)
\end{equation*}
exists, independently of the specific sequence of partitions and the choice of $\{{\boldsymbol \xi}_{\alpha}\}$, then $f$ is said to be Riemann integrable, and this limit is called the ($d$-tuple) Riemann integral of $f$, which we denote by 
$\displaystyle\int_{\mathbb{R}^d}f({\bf x})d{\bf x}$.

In particular, if we take the sequence of partitions given by$\Delta_{\epsilon}=\{\epsilon D_{\alpha}|~\alpha\in A\}$ ($\epsilon>0$), then, for every Riemann integrable function $f$, we have 
\begin{equation}\label{eq:riemannsumlimit1}
\lim_{\epsilon\to +0}\sum_{\alpha\in A}\epsilon^df(\epsilon {\boldsymbol \xi}_{\alpha}){\rm vol}(D_{\alpha})=\int_{\mathbb{R}^d}f({\bf x})d{\bf x},
\end{equation}
where we should note that ${\rm vol}(\epsilon D_{\alpha})=\epsilon^d{\rm vol}(D_{\alpha})$.

Now we look at Eq.~\!\ref{eq:riemannsumlimit1} from a different angle.
We think that $\omega(\boldsymbol{\xi}_{\alpha}):={\rm vol}(D_{\alpha})$ is a {\it weight}\index{weight} of the point $\boldsymbol{\xi}_{\alpha}$, and that Eq.~\!\ref{eq:riemannsumlimit1} is telling how the {\it weighted discrete set}\index{weighted discrete set} $(\Gamma,\omega)$ is distributed in $\mathbb{R}^d$; more specifically we may consider that Eq.~\!\ref{eq:riemannsumlimit1} implies {\it uniformity}\index{uniformity}, in a weak sense, of $(\Gamma,\omega)$ in $\mathbb{R}^d$. This view motivates us to propose the following definition.

In general, a weighted  discrete subset in $\mathbb{R}^d$ is a discrete set $\Gamma\subset \mathbb{R}^d$ with a function $\omega:\Gamma\rightarrow \mathbb{C}\backslash \{0\}$. 
Given a compactly supported function $f$ on $\mathbb{R}^d$, define the (generalized) {\it Riemann sum associated with $(\Gamma,\omega)$}\index{generalized Riemann sum} by setting
$$\displaystyle
\sigma(f,\Gamma,\omega)=\sum_{{\bf z}\in \Gamma}f({\bf z})\omega({\bf z}).
$$
In addition, we say that $(\Gamma,\omega)$ has {\it constant density}\index{constant density} $c(\Gamma,\omega)\neq 0$ (Marklof and Str\"{o}mbergsson \cite{mar}) if 
\begin{equation}
\displaystyle\lim_{\epsilon\to +0}\sigma(f^{\epsilon},\Gamma,\omega)\Big(=\lim_{\epsilon\to +0}\sum_{z\in \Gamma}\epsilon^df(\epsilon z)\omega(z)\Big)=
c(\Gamma,\omega)\int_{\mathbb{R}^d}f({\bf x})d{\bf x}
\end{equation}
holds for any bounded Riemann integrable function $f$ on $\mathbb{R}^d$ with compact support, where
$\displaystyle
f^{\epsilon}(x)=\epsilon^df(\epsilon x)
$; thus the weighted discrete set associated with a partition $\{D_{\alpha}\}$ and $\{\boldsymbol{\xi}_{\alpha}\}$ has constant density $1$.
In the case $\omega\equiv 1$, we write $\sigma(f,\Gamma)$ for $\sigma(f,\Gamma,\omega)$, and $c(\Gamma)$ for $c(\Gamma,\omega)$ when $\Gamma=(\Gamma, \omega)$ has constant density.

In connection with the notion of constant density, it is perhaps worth recalling the definition of a {\it Delone set}\index{Delone set}, a qualitative concept of ``uniformity". A discrete set $\Gamma$ is called a Delone set if it satisfies the following two conditions (Delone \cite{del}). 

\medskip

(1) There exists $R>0$ such that every ball $B_R(x)$ (of radius $R$ whose center is $x$) has a nonempty intersection with $\Gamma$, i.e., $\Gamma$ is {\it relatively dense}\index{relatively dense};

\smallskip

(2) there exists $r>0$ such that each ball $B_r(x)$ contains at most one element of $\Gamma$, i.e., $\Gamma$ is {\it uniformly discrete}\index{uniformly discrete}. 

\medskip

The following proposition states a relation between Delone sets and Riemann sums.

\begin{prop}\label{prop:delone} Let $\Gamma$ be a Delone set. Then there exist positive constants $c_1,c_2$ such that
$$
c_1\int_{\mathbb{R}^d}f({\bf x})d{\bf x}\leq
\varliminf_{\epsilon\to +0}\sigma(f^{\epsilon},\Gamma)
\leq \varlimsup_{\epsilon\to +0}\sigma(f^{\epsilon},\Gamma)\leq 
c_2\int_{\mathbb{R}^d}f({\bf x})d{\bf x}
$$
for every nonnegative-valued function $f$.
\end{prop}

\noindent{\it Proof}\quad In view of the Delone property, one can find two partitions $\{D_{\alpha}\}$ and $\{D'_{\beta}\}$ consisting of rectangular parallelotopes satisfying

\smallskip
(i)~ Every $D_{\alpha}$ has the same size, and contains at least one element of $\Gamma$;

\smallskip
(ii)~ every $D'_{\beta}$ has the same size, and contains at most one element of $\Gamma$.

\smallskip
Put $c_1={\rm vol}(D_{\alpha})^{-1}$ and $c_2={\rm vol}(D'_{\beta})^{-1}$. 
We take  a subset $\Gamma_1$ of $\Gamma$ such that every $D_{\alpha}$ contains just one element of $\Gamma_1$, and also take $\Gamma_2\supset \Gamma$ such that every $D'_{\beta}$ contains just one element of $\Gamma_2$. We then have
$\sigma(f^{\epsilon},\Gamma_1)\leq \sigma(f^{\epsilon},\Gamma)\leq \sigma(f^{\epsilon},\Gamma_2)$. Therefore using Eq.~\!\ref{eq:riemannsumlimit1}, we have
\begin{eqnarray*}
&&c_1\int_{\mathbb{R}^d}f({\bf x})d{\bf x}=\lim_{\epsilon\to +0}\sigma(f^{\epsilon},\Gamma_1)\leq \varliminf_{\epsilon\to +0}\sigma(f^{\epsilon},\Gamma)\\
&& \qquad \leq \varlimsup_{\epsilon\to +0}\sigma(f^{\epsilon},\Gamma)\leq \lim_{\epsilon\to +0}\sigma(f^{\epsilon},\Gamma_2)= c_2\int_{\mathbb{R}^d}f({\bf x})d{\bf x},
\end{eqnarray*}
where we should note that $\sigma_{\epsilon}(f,\Gamma_1)$ and $\sigma_{\epsilon}(f,\Gamma_2)$ are ordinary Riemann sums. 
\hfill $\Box$

\medskip
One might ask ``what is the significance of the notions of generalized Riemann sum and constant density?" Admittedly these notions are not so much profound (one can find more or less the same concepts in plural references). It may be, however, of great interest to focus our attention on the constant $c(\Gamma, \omega)$. In the subsequent sections, we shall give two ``arithmetical" examples for which the constant $c(\Gamma)$ is explicitly computed.

\section{Classical example 1}

Let $\mathbb{Z}_{\rm prim}^d$ $~(d\geq 2)$ be the set of {\it primitive lattice points}\index{primitive lattice point} in the $d$-dimensional standard lattice $\mathbb{Z}^d$, i.e., the set of {\it lattice points visible from the origin} (note that $\mathbb{Z}_{\rm prim}^2$ is the set of $(x,y)\in \mathbb{Z}^2$ such that $(|x|, |y|)$ is a coprime pair of positive integers, together with $(\pm 1,0)$ and $(0,\pm 1)$). 

\begin{thm}\label{thm:ppp}
$\mathbb{Z}_{\rm prim}^d$ has constant density $\zeta(d)^{-1}$; that is,
\begin{equation*}
\lim_{\epsilon\to +0}\sum_{{\bf z}\in \mathbb{Z}_{\rm prim}^d}\epsilon^df(\epsilon {\bf z})=\zeta(d)^{-1}\int_{\mathbb{R}^d}f({\bf x})d{\bf x}.
\end{equation*} 
Here $\displaystyle \zeta(s)=\sum_{n=1}^{\infty}n^{-s}$ is the {\it zeta function}\index{zeta function}. 
\end{thm}

The proof, which is more or less known as a sort of folklore, will be indicated in Sect.~\!\ref{sec:IEP}.

Noting that $\zeta(2)=\pi^2/6$ and applying this theorem to the indicator function $f$ for the square $\{(x,y)|~0\leq x,y\leq 1\}$, we obtain the following well-known statement.

\begin{cor}\label{gauss}
 The probability that two randomly chosen positive integers are coprime is $6/\pi^2$. More precisely
\begin{equation*} 
\lim_{N\to \infty}\frac{1}{N^2}\big|\big\{(a,b)\in \mathbb{N}\times \mathbb{N}|~{\rm gcd}(a,b)=1,~a,b\leq N\big\}\big|=\frac{6}{\pi^2},
\end{equation*}
where ${\rm gcd}(a,b)$ stands for the greatest common divisor of $a,b$.
\end{cor}

\vspace{-0.5cm}
\begin{figure}[htbp]
\begin{center}
\includegraphics[width=.55\linewidth]{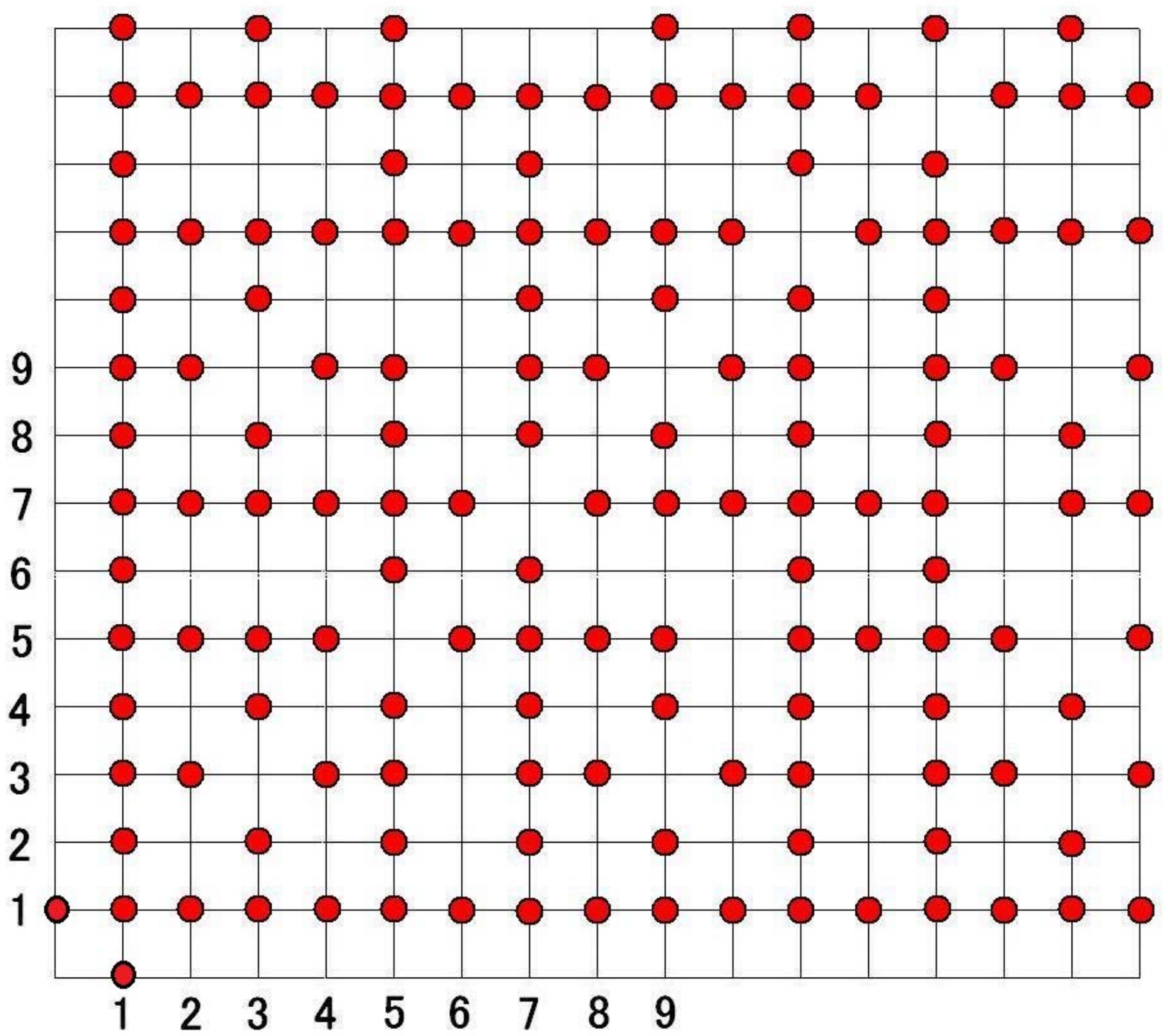}
\end{center}
\vspace{-4.3cm}
\caption{Coprime pairs}\label{fig:coprime}
\end{figure}

\begin{rem}{\rm 
(1)$~$ Gauss's {\it Mathematisches Tagebuch}\footnote{See vol. X in Gauss {\it Werke}.} (Mathematical Diary), a record of the mathematical discoveries of C. F. Gauss from 1796 to 1814, contains 146 entries, most of which consist of brief and somewhat cryptical statements. Some of the statements which he never published were independently discovered and published by others often many years later.\footnote{The first entry, the most famous one, records the discovery of the construction of a heptadecagon by ruler and compass. The diary was kept by Gauss's bereaved until 1899. It was St\"{a}ckel who became aware of the existence of the diary.}

The entry relevant to Corollary \ref{gauss} is the 31st dated 1796 September 6:

\medskip

``Numero fractionum inaequalium quorum denomonatores certum limitem non superant ad numerum fractionum omnium quarum num[eratores] 
aut denom[inatores] sint diversi infra limitem in infinito ut $6:\pi\pi$"

\medskip

This vague statement about counting (irreducible) fractions was formulated  in an appropriate way afterwards and proved rigorously by Dirichlet (1849) and Ernesto Ces\`{a}ro (1881). As a matter of fact, because of its vagueness, there are several ways to interpret what Gauss was going to convey.\footnote{For instance, see Ostwald's Klassiker der exakten Wissenschaften ; Nr. 256. The 14th entry dated 20 June, 1796 for which Dirichlet gave a proof is considered a companion of the 31st entry. The Yagloms \cite{yag} refer to the question on the probability of two random integers being coprime as ``Chebyshev's problem".}

\medskip
(2)$~$ In connection with Theorem \ref{thm:ppp}, it is perhaps worthwhile to make reference to the {\it Siegel mean value theorem}\index{Siegel mean value theorem} (\cite{siegel}).
 
Let $g\in \mathrm{SL}_d(\mathbb{R})$. For a bounded Riemann integrable function $f$ on $\mathbb{R}^d$ with compact support, we consider
$$
\varPhi(g)=\sum_{{\bf z}\in \mathbb{Z}^d\backslash\{{\bf 0}\}}f(g{\bf z}),\quad
\Psi(g)=\sum_{{\bf z}\in \mathbb{Z}_{{\rm prim}}^d}f(g{\bf z}).
$$
Both functions $\varPhi$ and $\Psi$ are $\mathrm{SL}_d(\mathbb{Z})$-invariant with respect to the 
right action of $\mathrm{SL}_d(\mathbb{Z})$ on $\mathrm{SL}_d(\mathbb{R})$, so that these are identified with functions on the coset space $\mathrm{SL}_d(\mathbb{R})/\mathrm{SL}_d(\mathbb{Z})$. Recall that $\mathrm{SL}_d(\mathbb{R})/\mathrm{SL}_d(\mathbb{Z})$ has finite volume with respect to the measure $dg$ induced from the Haar measure on $\mathrm{SL}_d(\mathbb{R})$. We assume $\displaystyle\int_{\mathrm{SL}_d(\mathbb{R})/\mathrm{SL}_d(\mathbb{Z})}1~dg=1$. Then the Siegel theorem asserts 
\begin{eqnarray*}
&&\int_{\mathrm{SL}_d(\mathbb{R})/\mathrm{SL}_d(\mathbb{Z})}\Big(\sum_{{\bf z}\in \mathbb{Z}^d\backslash\{{\bf 0}\}}f(g{\bf z})\Big)dg=\int_{\mathbb{R}^d}f({\bf x})d{\bf x},\label{eq:siegel}\\
&&\int_{\mathrm{SL}_d(\mathbb{R})/\mathrm{SL}_d(\mathbb{Z})}\Big(\sum_{{\bf z}\in \mathbb{Z}_{{\rm prim}}^d}f(g{\bf z})\Big)dg=\zeta(d)^{-1}\int_{\mathbb{R}^d}f({\bf x})d{\bf x}.\nonumber
\end{eqnarray*}
\hfill$\Box$

}
\end{rem}

\section{Classical example 2}

A {\it Pythagorean triple}\index{Pythagorean triple},\footnote{Pythagorean triples have a long history since the Old Babylonian period in Mesopotamia nearly 4000 years ago. Indeed, one can read 15 Pythagorean triples in the ancient tablet, written about 1800 BCE, called Plimpton 322 (Weil \cite{weil}).} the name stemming from the Pythagorean theorem for right triangles,  is a triple of positive integers $(\ell,m,n)$ satisfying the equation $\ell^2+m^2=n^2$. Since $(\ell/n)^2+(m/n)^2=1$, a Pythagorean triple yields a {\it rational point}\index{rational point} $(\ell/n,m/n)$ on the unit circle $S^1=\{(x,y)|~x^2+y^2=1\}$. Conversely any rational point on $S^1$ is derived from a Pythagorean triple. Furthermore the well-known parameterization of $S^1$ given by $x=(1-t^2)/(1+t^2)$, $y=2t/(1+t^2)$ tells us that the set of rational points $S^1(\mathbb{Q})=S^1\cap \mathbb{Q}^2$ is dense in $S^1$ (we shall see later how rational points are distributed from a quantitative viewpoint).

A Pythagorean triple $(x,y,z)$ is called {\it primitive}\index{primitive Pythagorean triple} if  $x,y,z$  are pairwise coprime. ``Primitive" is so named because any Pythagorean triple is generated trivially from the primitive one, i.e., if $(x,y,z)$ is Pythagorean, there are a positive integer $\ell$ and a primitive $(x_0,y_0,z_0)$ such that $(x,y,z)=(\ell x_0,\ell y_0,\ell z_0)$.

The way to produce primitive Pythagorean triples (PPTs) is described as follows: 
If $(x,y,z)$ is a PPT, then there exist positive integers $m,n$ such that

\smallskip

{\rm (i)}~ $m>n$, 

\smallskip

{\rm (ii)}~ $m$ and $n$ are coprime,

\smallskip

{\rm (iii)}~ $m$ and $n$ have different parity,

\smallskip

{\rm (iv)} $~(x,y,z)=(m^2-n^2,2mn,m^2+n^2)$ or
$(x,y,z)$ $=$ $(2mn,m^2-n^2,m^2+n^2)$. 

\medskip
Conversely, if $m$ and $n$ satisfy {\rm (i), (ii), (iii)}, then $(m^2-n^2,2mn,m^2+n^2)$ and $(2mn,m^2-n^2,m^2+n^2)$ are PPTs.

In the table below, due to M. Somos \cite{somos}, of PPTs $(x,y,z)$ enumerated in ascending order with respect to $z$, the triple $(x_N,y_N,z_N)$ is the $N$-th PPT (we do not discriminate between $(x,y,z)$ and $(y,x,z)$).

\begin{table}[hbtp]
\begin{tabular}{|c||c|c|c|}\hline
$N$ & $x_N$ & $y_N$ & $z_N$ \\ \hline  
1   &   3  &    4  &    5 \\    

2   &   5  &   12  &   13 \\   

3   &  15 &     8  &   17 \\  

4   &   7 &    24  &   25 \\   

5   &  21 &    20 &    29 \\   

6   &  35  &   12 &    37 \\   

7   &   9  &   40 &    41 \\   

8   &  45  &   28  &   53 \\    

9   &  11  &   60  &   61 \\   

10  &   63  &   16 &    65 \\
\hline
\end{tabular}
\hspace{0.1cm}
\begin{tabular}{|c||c|c|c|}\hline
$N$ & $x_N$ & $y_N$ & $z_N$ \\ \hline 
 11 &    33 &    56 &    65\\   
 12 &    55 &    48 &    73\\   
 13 &    77 &    36 &    85\\   
 14 &    13 &    84 &    85\\   
 15 &    39 &    80 &    89\\   
 16 &    65 &    72 &    97\\  
 17 &    99 &    20 &   101\\  
 18 &    91 &    60 &   109\\    
 19 &    15 &   112 &   113\\    
 20 &   117 &    44 &   125\\
\hline
\end{tabular} $\cdots$ 
\begin{tabular}{|c||c|c|c|}\hline
$N$ & $x_N$ & $y_N$ & $z_N$ \\ \hline  
1491  & 4389  & 8300  & 9389\\   

1492  &  411  & 9380  & 9389\\   

1493  &  685  & 9372  & 9397\\   

1494  &  959   & 9360  & 9409\\   

1495  & 9405  &  388  & 9413\\    

1496  & 5371  & 7740  & 9421\\   

1497  & 9393  &  776  & 9425\\    

1498  & 7503  & 5704  & 9425\\   

1499  & 6063  & 7216 &  9425\\   

1500  & 1233  & 9344  & 9425\\   
\hline
\end{tabular}
\end{table}

What we have interest in is the {\it asymptotic behavior} of $z_N$ as $N$ goes to infinity. The numerical observation tells us that the sequence $\{z_N\}$ almost linearly increases as $N$ increases. Indeed
$
z_{100}/100=6.29,~z_{1000}/1000=6,277,
z_{1500}/1500=6.28333\cdots$,
which convinces us that $\displaystyle\lim_{N\to \infty}z_N/N$ exists (though the speed of convergence is very slow), and the limit is expected to be equal to $2\pi=6.2831853\cdots$.  
This is actually true as shown by D. N. Lehmer \cite{lehmer} in 1900, though his proof is by no means easy. 

We shall prove Lehmer's theorem\index{Lehmer's theorem} by counting coprime pairs $(m,n)$ satisfying the condition that $m-n$ is odd. A key of our proof is the following theorem.

\begin{thm}\label{thm:ppt1} $\mathbb{Z}_{\rm prim}^{2,*}=\{(m.n)\in \mathbb{Z}_{\rm prim}^2|~m-n~\text{is odd}\}~(=\{(m.n)\in \mathbb{Z}_{\rm prim}^2|~m-n\equiv1~\text{(\rm mod~2})\})$ has constant density $4/\pi^2$; namely
\begin{equation}\label{eq:limit2}
\lim_{\epsilon\to +0}\sum_{{\bf z}\in \mathbb{Z}_{\rm prim}^{2,*}}\epsilon^2f(\epsilon {\bf z})=\frac{2}{3}\zeta(2)^{-1}\int_{\mathbb{R}^2}f({\bf x})d{\bf x}=\frac{4}{\pi^2}\int_{\mathbb{R}^2}f({\bf x})d{\bf x}.
\end{equation}

\end{thm}

We postpone the proof to Sect.~\!\ref{sec:IEP}, and apply this theorem to the indicator function $f$ for the set $\{(x,y)\in \mathbb{R}^2|~\!x\geq y,~x^2+y^2\leq 1\}$. Since
\begin{eqnarray*}
&&\sum_{{\bf z}\in \mathbb{Z}_{\rm prim}^{2,*}}\epsilon^2f(\epsilon {\bf z})=\epsilon^2
\big|
\big\{(m,n)\in \mathbb{N}^2|~
 {\rm gcd}(m,n)=1,~m>n,\\
&& \qquad\qquad \qquad \qquad m^2+n^2\leq \epsilon^{-2},
~m-n~\text{is odd}\big\}
\big|,
\end{eqnarray*}
we obtain
\begin{eqnarray}\label{eq:limitx3}
&&\lim_{N\to\infty}\frac{1}{N}\big|\big\{(m,n)\in \mathbb{N}^2|~
 {\rm gcd}(m,n)=1,~m>n,~m^2+n^2\leq N,\\
&&\qquad \qquad \qquad m-n~\text{is odd}\big\}\big|= \frac{2}{3}\cdot\frac{6}{\pi^2}\cdot\frac{\pi}{8}=\frac{1}{2\pi}.\nonumber
\end{eqnarray}
Note that $\big|\big\{(m,n)\in \mathbb{N}^2|~{\rm gcd}(m,n)=1,~m>n,~m^2+n^2\leq N,~ m-n~\text{is odd}\big\}\big|$ coincides with the number of PPT $(x,y,z)$ with $z\leq N$. This observation leads us to 

\begin{cor} {\rm (Lehmer)}\quad 
$\displaystyle\lim_{N\to\infty} \frac{z_N}{N}=2\pi$. 
\end{cor}

\begin{rem}{\rm  
Fermat's theorem on sums of two squares,\footnote{Every prime number $p =
4k + 1$ is in one and only one way a sum of two squares of positive integers.} together with his {\it little theorem} and the formula $(a^2+b^2)(c^2+d^2)=(ac\pm bd)^2+(ad\mp bc)^2$, yields the following complete characterization of PPTs which is substantially equivalent to the result stated in the letter from  Fermat to Mersenne dated 25 December 1640 (cf.~\!Weil \cite{weil}). 

\medskip
{\it 
An odd number $z$ is written as $m^2+n^2$ by using two coprime positive integers $m, n$ (thus automatically having different parity) if and only if every prime divisor of $z$ is of the form $4k+1$. In other words, the set 
$\{z_N\}$ coincides with the set of odd numbers whose prime divisors are of the form $4k+1$. Moreover, if we denote by $\nu(z)$ the number of distinct prime divisors of $z$, then
$z=z_N$ in the list is repeated $2^{\nu(z)-1}$ times.}
\hfill$\Box$

}
\end{rem}

Theorem \ref{thm:ppt1} can be used to establish

\begin{cor}\label{cor:uni}
\quad 
For a rational point $(p,q)\in S^1(\mathbb{Q}) (=S^1\cap\mathbb{Q}^2)$, define the {\it height} $h(p,q)$ to be the minimal positive integer $h$ such that $(hp,hq)\in \mathbb{Z}^2$. Then for any arc $A$ in $S^1$, we have
$$
\big|\big\{(p,q)\in A\cap\mathbb{Q}^2|~h(p,q)\leq h\big\}\big|\sim \frac{2\cdot{\rm length}(A)}{\pi^2}h\quad (h\to \infty),
$$
and hence rational points are {\it equidistributed}\index{equidistributed} on the unit circle in the sense that 
$$
\lim_{h\to \infty}\frac{\big|\big\{(p,q)\in A\cap\mathbb{Q}^2|~h(p,q)\leq h\big\}\big|}{\big|\big\{(p,q)\in S^1\cap\mathbb{Q}^2|~h(p,q)\leq h\}\big|}=\frac{{\rm length}(A)}{2\pi}.
$$

\end{cor}

In his paper \cite{duke}, W. Duke suggested that this corollary can be proved by using tools from the theory of $L$-functions combined with Weyl's famous criterion for {\it equidistribution on the circle} (\cite{weyl}). Our proof below relies on a generalization of Eq.~\!\ref{eq:limitx3}.

Given $\alpha, \beta$ with $0\leq \alpha<\beta\leq 1$, we put 
\begin{eqnarray*}
&&P(N;\alpha,\beta)=\big\{(m,n)\in \mathbb{N}\times\mathbb{N}|~{\rm gcd}(m,n)=1, \alpha\leq n/m\leq \beta,\\
&& \qquad \qquad m-n~\text{is odd},~ m^2+n^2\leq N\big\}.
\end{eqnarray*}
Namely we count coprime pairs $(m,n)$ with odd $m-n$ in the circular sector 
$$
\{(x,y)\in \mathbb{R}^2|~ x, y>0,~\alpha x\leq y\leq \beta x,~x^2+y^2\leq N\}.
$$
Since the area of the region $
\{(x,y)\in \mathbb{R}^2|~ x, y>0,~\alpha x\leq y\leq \beta x,~x^2+y^2\leq 1\}
$ is $\displaystyle\frac{1}{2}\arctan \frac{\beta-\alpha}{1+\alpha\beta}$, applying again Eq.~\!\ref{eq:limit2} to the indicator function for this region, we obtain$$
\displaystyle\lim_{N\to\infty}\frac{1}{N}\big|P(N;\alpha,\beta)\big|= \frac{2}{\pi^2} \arctan \frac{\beta-\alpha}{1+\alpha\beta}.
$$

Now we sort points 
$(p,q)\neq (\pm 1,0), (0,$ $\pm 1)$ in $S^1(\mathbb{Q})$ by 4 quadrants containing $(p,q)$, and also by parity of $x$ when we write $|p|=x/z$, $|q|=y/z$ with a PPT $(x,y,z)$. Here we should notice that $h(p,q)=z=m^2+n^2$. 
Thus counting rational points with the height function $h(p,q)$ reduces to counting PPTs. 

Put 
\begin{eqnarray*}
&&S_{\mathbb{Q}}^1({\rm odd})=\big\{(p,q)\in S^1(\mathbb{Q})|~\text{$x$ is odd}\big\},\\
&&S_{\mathbb{Q}}^1({\rm even})=\big\{(p,q)\in S^1(\mathbb{Q})|~\text{$x$ is even}\big\}.
\end{eqnarray*}
Then
$$
S^1(\mathbb{Q})=S_{\mathbb{Q}}^1({\rm odd})\cup S_{\mathbb{Q}}^1({\rm even})\cup \big\{(
\pm 1, 0), (0, \pm 1)
\big\}\quad (\text{disjoint}).
$$
Note that the correspondence $(p,q)\mapsto (q,p)$ interchanges $S_{\mathbb{Q}}^1({\rm odd})$ and $S_{\mathbb{Q}}^1({\rm even})$. 
Therefore, in order to complete the proof, it is enough to show that 
\begin{eqnarray*}
&&\big|\big\{(p,q)\in S_{\mathbb{Q}}^1({\rm odd})|~\theta_1\leq \theta(p,q)<\theta_2,~h(p,q)\leq h\big\}\big|\\
&\sim& \frac{1}{\pi^2}(\theta_2-\theta_1)h \quad (h\to\infty),\nonumber
\end{eqnarray*}
where $(p,q)=\big(\cos \theta(p,q),\sin \theta(p,q)\big)$.  
Without loss of generality, one may assume $0\leq \theta_1<\theta_2\leq \pi/2$. Since
$$
\tan \theta(p,q)=\frac{q}{p}=\frac{2mn}{m^2-n^2}=\frac{2\displaystyle\frac{n}{m}}{1-\Big(\displaystyle\frac{n}{m}\Big)^2},
$$
if we define $\Theta(m,n)\in [0,\pi/2)$ by $\tan \Theta(m,n)=n/m$, then 
$
\tan \theta(p,q)= \tan 2\Theta(m,n)
$,
and hence $\theta(p,q)=2\Theta(m,n)$. Therefore 
\begin{eqnarray*}
&&\big|\big\{(p,q)\in S_{\mathbb{Q}}^1({\rm odd})|~\theta_1\leq \theta(g)<\theta_2,~h(p,q)\leq h\big\}\big|\\
&=&\big|P(h; \arctan \theta_1/2, \arctan \theta_2/2)\big|\\
&\sim & \frac{1}{\pi^2}(\theta_2-\theta_1)h,
\end{eqnarray*}
as required.

\begin{rem}
{\rm 
Interestingly, $S^1(\mathbb{Q})$ (and hence Pythagorean triples)  has something to do with crystallography. 
Indeed $S^1(\mathbb{Q})$ with the natural group operation 
is an example of {\it coincidence symmetry groups} that show up in the theory of crystalline interfaces and grain boundaries\footnote{Grain boundaries are interfaces where crystals of different orientations meet.} in polycrystalline materials (Ranganathan \cite{rang}, Zeiner \cite{zeiner}). This theory is concerned with partial coincidence of lattice points in two identical crystal lattices. See \cite{sunada} for the details, and also \cite{sunada1} for the mathematical theory of crystal structures.
\hfill$\Box$

}
\end{rem}

\section{The Inclusion-Exclusion Principle}\label{sec:IEP}
The proof that the discrete sets $\mathbb{Z}^d_{\rm prim}$ and $\mathbb{Z}^{2,*}_{\rm prim}$ have constant density relies on the identities derived from  the so-called {\it Inclusion-Exclusion Principle}\index{Inclusion-Exclusion Principle} (IEP), which is a generalization of the obvious equality $|A\cup B|=|A|+|B|-|A\cap B|$  for two finite sets $A, B$. Despite its simplicity, the IEP is a powerful tool to approach general {\it counting problems} involving aggregation of things that are not mutually exclusive (Comtet \cite{com}).

To state the IEP in full generality, we consider a family $\{A_h\}_{h=1}^{\infty}$ of subsets of $X$ where $X$ and $A_h$ are not necessarily finite. 
Let $f$ be a real-valued function with finite support defined on $X$. We assume that there exists $N$ such that if $h>N$, then $A_h\cap {\rm supp}~\!f=\emptyset$, i.e. $f(x)=0$ for $x\in A_h$. In the following theorem, the symbol $A^c$ means the complement of a subset $A$ in $X$.

\begin{thm}\label{thm:thm}{\rm (}Inclusion-Exclusion Principle{\rm )}
\begin{eqnarray}\label{eq:exclu}
\sum_{x\in \bigcap_{h=1}^{\infty}A_h^c}f(x)&=&\sum_{k=0}^{\infty}(-1)^k\sum_{h_1<\cdots<h_k}
\sum_{x\in A_{h_1}\cap\cdots \cap A_{h_k}}f(x)\\
&&\left(=\sum_{k=0}^{N}(-1)^k\sum_{h_1<\cdots<h_k}
\sum_{x\in A_{h_1}\cap\cdots \cap A_{h_k}}f(x)\right),\nonumber
\end{eqnarray}
where, for $k=0$, the term $\displaystyle\sum_{h_1<\cdots<h_k}
\sum_{x\in A_{h_1}\cap\cdots \cap A_{h_k}}f(x)$ should be understood to be 
$\displaystyle\sum_{x\in X}f(x)$.

\end{thm}

For the proof, one may assume, without loss of generality, that $X$ is finite, and it suffices to handle the case of a finite family $\{A_h\}_{h=1}^N$. The proof is accomplished by induction on $N$.

Making use of the IEP, we obtain the following theorem (this is actually an easy exercise of the IEP; see Vinogradov \cite{vino} for instance). 

\begin{thm}\label{eq:vino}
\begin{equation*}\label{eq:vinox}
\sum_{{\bf z}\in \mathbb{Z}^d_{\rm prim}}f({\bf z})=\sum_{k=1}^{\infty}\mu(k)\sum_{{\bf w}\in \mathbb{Z}^d\backslash\{{\bf 0}\}}f(k{\bf w}),
\end{equation*}
where  $f$ is a function on $\mathbb{R}^d$ with compact support (thus both sides are finite sums), and 
$\mu(k)$ is the {\it M\"{o}bius function}:
$$
\mu(k)=\begin{cases}
1& (k=1)\\
(-1)^r & (k=p_{h_1}\cdots p_{h_r};~h_1<\cdots<h_r)\\
0& (\text{otherwise}),
\end{cases}
$$
where $p_1<p_2<\cdots$ are all primes enumerated into ascending order.
\end{thm}

The proof goes as follows. Consider the case that  
$$
X=\mathbb{Z}^d\backslash\{{\bf 0}\},\quad A_h=\{(x_1,\ldots, x_d)\in X|~p_h|x_1,\ldots, p_h|x_d\}.
$$
Then $\displaystyle\bigcap_{h=1}^{\infty} A_h^{c}=\mathbb{Z}_{\rm prime}^d$. We also easily observe 
$$
A_{h_1}\cap\cdots \cap A_{h_k}=p_{h_1}\cdots p_{h_k}X.
$$
Applying Eq.~\!\ref{eq:exclu} to this case, we have
\begin{eqnarray*}
\sum_{{\bf z}\in\mathbb{Z}^d_{\rm prime}}f({\bf z})&=&\sum_{k=0}^{\infty}(-1)^k\sum_{h_1<\cdots<h_k}\sum_{{\bf w}\in \mathbb{Z}^d\backslash\{{\bf 0}\}}f(p_{h_1}\cdots p_{h_k}{\bf w})\\
&=& \sum_{k=1}^{\infty}\mu(k)\sum_{{\bf w}\in \mathbb{Z}^d\backslash\{{\bf 0}\}}f(k{\bf w}).
\end{eqnarray*}

\noindent {\bf Proof of Theorem \ref{thm:ppp}} Applying Theorem \ref{eq:vino} to $f^{\epsilon}$, we have
$$
\sum_{{\bf z}\in \mathbb{Z}^d_{\rm prim}}\epsilon^df(\epsilon{\bf z})=\sum_{k=1}^{\infty}\mu(k)\sum_{{\bf w}\in \mathbb{Z}^d\backslash\{{\bf 0}\}}\epsilon^df(\epsilon k{\bf w}),
$$
What we have to confirm is the exchangeability of the limit 
and summation:
$$
\lim_{\epsilon\to +0}\sum_{k=1}^{\infty}\left(\mu(k)\sum_{{\bf w}\in \mathbb{Z}^d\backslash\{{\bf 0}\}}\epsilon^df(\epsilon k{\bf w})\right)=\sum_{k=1}^{\infty}\lim_{\epsilon\to +0}\left(\mu(k)\sum_{{\bf w}\in \mathbb{Z}^d\backslash\{{\bf 0}\}}\epsilon^df(\epsilon k{\bf w})\right).
$$  
If we take this for granted, then we easily get the claim since 
$$
\lim_{\epsilon\to +0}\sum_{{\bf w}\in \mathbb{Z}^2\backslash\{{\bf 0}\}}\epsilon^df(\epsilon k{\bf w})=k^{-d}\lim_{\delta\to +0}\sum_{{\bf w}\in \mathbb{Z}^d\backslash\{{\bf 0}\}}\delta^df(\delta{\bf w})=k^{-d}\int_{\mathbb{R}^d}f({\bf x})d{\bf x},
$$
and
$\displaystyle
\sum_{k=1}^{\infty}\mu(k)k^{-d}=\zeta(d)^{-1}
$. As a matter of fact, the exchangeability does not follow from Weierstrass' M-test in a direct manner.  One can check it by a careful argument.\hfill$\Box$

\medskip

In the case of Theorem \ref{thm:ppt1}, we consider
$$
\big(\mathbb{Z}^{\rm odd})^2_{\rm prim}=
\big\{(m,n)\in \mathbb{Z}^{\rm odd}\times \mathbb{Z}^{\rm odd}|~{\rm gcd}(m,n)=1\big\},
$$
where $\mathbb{Z}^{\rm odd}$ is the set of odd integers. 
Then 
$$
\mathbb{Z}_{\rm prim}^{2,*}=\mathbb{Z}_{\rm prim}^2\backslash \big(\mathbb{Z}^{\rm odd})^2_{\rm prim}.
$$
Therefore it suffices to show that $\big(\mathbb{Z}^{\rm odd})^2_{\rm prim}$ has constant density $2/\pi^2$. This is done by using the following theorem for which we need a slightly sophisticated use of the IEP. 

\begin{thm}\label{thm:pptx1}
\begin{equation*}\label{eq:pptx}
\sum_{{\bf z}\in (\mathbb{Z}^{\rm odd})^2_{\rm prim}}f({\bf z})=
\sum_{k=1}^{\infty}\mu(k)\sum_{h=0}^{\infty}\sum_{{\bf w}\in (\mathbb{Z}^{\rm odd})^2}
f(k 2^h{\bf w}).
\end{equation*}
\end{thm}

For the proof, we put
$$
X=\coprod_{\ell=1}^{\infty}\ell(\mathbb{Z}^{\rm odd})^2_{\rm prim}, \quad
A_h=\big\{(x,y)\in X|~ p_h|x~\text{and}~p_h|y\big\}.
$$

\begin{lemma}
$\displaystyle
A_{h_1}\cap\cdots\cap A_{h_k}=\coprod_{h=0}^{\infty}p_{h_1}\cdots p_{h_k}2^{h}(\mathbb{Z}^{\rm odd})^2
$.

\end{lemma}

\noindent{\it Proof}\quad It suffices to prove that $\displaystyle
A_{h_1}\cap\cdots\cap A_{h_k}=p_{h_1}\cdots p_{h_k}X
$ since any positive integer $\ell$ is expressed as $2^i\times {\rm odd}$. Clearly $A_{h_1}\cap\cdots\cap A_{h_k}\supset p_{h_1}\cdots p_{h_k}X$. Let $(x,y)\in A_{h_1}\cap\cdots\cap A_{h_k}$. Then one can find $(a,b)\in \mathbb{Z}^2$ such that $x=p_{h_1}\cdots p_{h_k}a$ and $y=p_{h_1}\cdots p_{h_k}b$. Moreover there exist $\ell\in \mathbb{N}$ and $(m,n)\in (\mathbb{Z}^{\rm odd})^2_{\rm prim}$ such that $x=\ell m$, $y=\ell n$, so 
$p_{h_1}\cdots p_{h_k}|{\rm gcd}(\ell m,\ell n)=\ell$. Therefore
$(x,y)\in p_{h_1}\cdots p_{h_k}X$.\hfill$\Box$

\begin{lemma}
$\displaystyle
\left(\bigcup_{h=1}^{\infty}A_h\right)^c=(\mathbb{Z}^{\rm odd})^2_{\rm prim}.
$
\end{lemma}

\noindent{\it Proof}\quad Obviously $\displaystyle\coprod_{\ell=2}^{\infty} \ell(\mathbb{Z}^{\rm odd})^2_{\rm prim}= \bigcup_{h=1}^{\infty}A_h$, from which the claim follows.\hfill$\Box$

\medskip

Theorem \ref{thm:pptx1} is a consequence of the above two lemmas.

Now using Theorem \ref{thm:pptx1}, we have
\begin{eqnarray*}
\sum_{{\bf z}\in (\mathbb{Z}^{\rm odd})^2_{\rm prim}}\epsilon^2f(\epsilon{\bf z})&=&
\sum_{k=0}^{\infty}(-1)^k\sum_{h_1<\cdots<h_k}\sum_{h=0}^{\infty}\sum_{{\bf w}\in (\mathbb{Z}^{\rm odd})^2}
\epsilon^2f(\epsilon p_{h_1}\cdots p_{h_k}2^h{\bf w})\\
&=& \sum_{k=1}^{\infty}\mu(k)\sum_{h=0}^{\infty}\sum_{{\bf w}\in (\mathbb{Z}^{\rm odd})^2}
\epsilon^2f(\epsilon k 2^h{\bf w}).
\end{eqnarray*}
We also have 
$$
\lim_{\epsilon\to +0}\sum_{{\bf z}\in (\mathbb{Z}^{\rm odd})^2}\epsilon^2f(\epsilon{\bf z})=\frac{1}{4}\int_{\mathbb{R}^2}f({\bf x})d{\bf x},
$$
since the left-hand side is the ordinary Riemann sum associated with the partition by the squares with side length $2$, 
and hence 
$$
\lim_{\epsilon\to +0}
\sum_{{\bf w}\in (\mathbb{Z}^{\rm odd})^2}\epsilon^2f(\epsilon k2^h{\bf w})=\frac{1}{(k2^h)^2}\frac{1}{4}\int_{\mathbb{R}^2}f({\bf x})d{\bf x}.
$$
Thus 
\begin{eqnarray*}
&&\lim_{\epsilon\to +0}\sum_{{\bf z}\in (\mathbb{Z}^{\rm odd})^2_{\rm prim}}\epsilon^2f(\epsilon {\bf z})\\
&=&\zeta(2)^{-1}\sum_{h=0}^{\infty}\frac{1}{4^h}\int_{\mathbb{R}^2}f({\bf x})d{\bf x}=\frac{6}{\pi^2}\cdot\frac{1}{3}\int_{\mathbb{R}^2}f({\bf x})d{\bf x}=\frac{2}{\pi^2}\int_{\mathbb{R}^2}f({\bf x})d{\bf x},
\end{eqnarray*}
as desired (this time, the exchangeability of the limit 
and summation is confirmed by Weierstrass' M-test). \hfill$\Box$

\begin{rem}{\rm 
Historically IEP was, for the first time, employed by Nicholas Bernoulli (1687--1759) to solve a combinatorial problem related to permutations.\footnote{The probabilistic form of IEP is attributed to de Moivre (1718). Sometimes IEP is referred to as the formula of Da Silva, or Sylvester.} More specifically he counted the number of {\it derangements}, that is, permutations such that none of the elements appears in its original position.\footnote{This problem (``probl\`{e}me des rencontres") was proposed by Pierre Raymond de Montmort in 1708. He solved it in 1713 at about the same time as did N. Bernoulli.} His result is pleasingly phrased, in a similar fashion as in the case of coprime pairs, as ``the probability that randomly chosen permutations are derangements is $1/e$" ($e$ is the base of natural logarithms).\hfill$\Box$

}
\end{rem}

\section{Generalized Poisson summation Formulas}
Generalized Riemann sums appear in the theory of {\it quasicrystals}\index{quasicrystals}, a form of solid matter whose atoms are arranged like those of a crystal but assume patterns that do not exactly repeat themselves. 

The interest in quasicrystals arose when in 1984 Schechtman et al.~\!\cite{SHB} discovered materials whose X-ray diffraction spectra had sharp spots indicative of long range order. Soon after the announcement of their discovery, material scientists and mathematicians began intensive studies of quasicrystals from both the empirical and theoretical sides.\footnote{As will be explained below, the theoretical discovery of quasicrystal structures was already made by R. Penrose in 1973.  
See Senechal and Taylor \cite{sene} for an account on the theory of quasicrystals at the early stage.} 

At the moment, there are several ways to mathematically define quasicrystals (see Lagarias \cite{lag} for instance). As a matter of fact,  an official nomenclature has not yet been agreed upon. 
In many reference, however, the Delone property for the discrete set $\Gamma$ representing the location of atoms is adopted as a minimum requirement for the characterization of quasicrystals. In addition to the Delone property, many authors assume that a {\it generalized Poisson summation formula}\index{generalized Poisson summation formula} holds for $\Gamma$, which embodies the patterns of X-ray diffractions for a real quasicrystal.

Let us recall the classical Poisson summation formula\index{Poisson summation formula}. For a {\it lattice group} $L$, a subgroup of $\mathbb{R}^d$ generated by a basis of $\mathbb{R}^d$, we denote by 
$L^{*}$ the dual lattice of $L$, i.e., $L^*=\{\boldsymbol{\eta}\in \mathbb{R}^d|~\langle \boldsymbol{\eta},{\bf z}\rangle \in \mathbb{Z}~\text{for every}~{\bf z}\in L\}$, and also denote by $D_{L}$ a fundamental domain for $L$.  We then have 
\begin{equation}\label{eq:poisson}
\sum_{{\bf z}\in L}f({\bf z})e^{2\pi i\langle {\bf z},\boldsymbol{\eta}\rangle}
={\rm vol}(D_{L})^{-1}\sum_{\boldsymbol{\xi}\in L^*}\hat{f}(\boldsymbol{\xi}-\boldsymbol{\eta}) \quad (i=\sqrt{-1}),
\end{equation}
in particular, 
\begin{equation}\label{eq:poisson1}
\sum_{{\bf z}\in L}f({\bf z})
={\rm vol}(D_{L})^{-1}\sum_{\boldsymbol{\xi}\in L^*}\hat{f}(\boldsymbol{\xi}),
\end{equation}
which is what we usually call the Poisson summation formula.  
Here $\hat{f}$ is the {\it Fourier transform} of a rapidly decreasing smooth function $f$:
$$
\hat{f}(\boldsymbol{\xi})=\int_{\mathbb{R}^d}f({\bf x})e^{-2\pi i\langle {\bf x},\boldsymbol{\xi}\rangle}d{\bf x}.
$$
Note that the left-hand side of Eq.~\!\ref{eq:poisson} is the Riemann sum $\sigma(f,L,\omega_{\boldsymbol{\eta}})$ for the weighted discrete set $(L,\omega_{\boldsymbol{\eta}})$, where $\omega_{\boldsymbol{\eta}}({\bf z})=e^{2\pi i\langle {\bf z},\boldsymbol{\eta}\rangle}$.

Having Eq.~\!\ref{eq:poisson1} in mind, we say that a generalized Poisson formula holds for $\Gamma$ if there exist a countable subset $\Lambda\subset \mathbb{R}^d$ and a sequence $\{a(\boldsymbol{\xi})\}_{\boldsymbol{\xi}\in \Lambda}$ such that 
\begin{equation}\label{eq:poisson2}
\sum_{{\bf z}\in \Gamma}f({\bf z})=\sum_{\boldsymbol{\xi}\in \Lambda}a(\boldsymbol{\xi})\hat{f}(\boldsymbol{\xi})
\end{equation}
for every compactly supported smooth function $f$. 

What we must be careful about here is that the set $\Lambda$ is allowed to have accumulation points, so that one cannot claim that the right-hand side of Eq.~\!\ref{eq:poisson2} converges in the ordinary sense. Thus the definition above is rather formal. One of the possible justifications is to assume that there exist an increasing family of subsets $\{\Lambda_N\}_{N=1}^{\infty}$ and functions 
$a_N(\boldsymbol{\xi})$ defined on $\Lambda_N$ such that

\medskip

(i)\quad $\displaystyle\bigcup_{N=1}^{\infty}\Lambda_N=\Lambda$,
  
\smallskip
(ii)\quad $\displaystyle\sum_{\boldsymbol{\xi}\in \Lambda_N}a_N(\boldsymbol{\xi})\hat{f}(\boldsymbol{\xi})$ converges absolutely,

\smallskip

(iii)\quad $\displaystyle\lim_{N\to\infty}a_N(\boldsymbol{\xi})=a(\boldsymbol{\xi})$,

\smallskip
(iv)\quad  $\displaystyle\sum_{{\bf z}\in \Gamma}f({\bf z})=\lim_{N\to\infty}\sum_{\boldsymbol{\xi}\in \Lambda_N}a_N(\boldsymbol{\xi})\hat{f}(\boldsymbol{\xi})$.

\medskip
We shall say that a discrete set $\Gamma$ is a {\it quasicrystal of Poisson type}\index{quasicrystal of Poisson type} if a generalized Poisson formula holds for $\Gamma$.\footnote{Some people use the term ``Poisson comb" in a bit different formulation.}

A typical class of quasicrystals of Poisson type is constructed by the {\it cut and project method}\index{cut and project method}.\footnote{This method was invented by de Bruijn \cite{debr}, and developed by many authors.} 
Let $L$ be a lattice group in $\mathbb{R}^N=\mathbb{R}^d\times \mathbb{R}^{N-d}$ ($N>d$), and let $W$ be a compact domain (called a {\it window}) in $\mathbb{R}^{N-d}$. We denote by $p_d$ and $p_{N-d}$ the orthogonal projections of $\mathbb{R}^N$ onto $\mathbb{R}^d$ and $\mathbb{R}^{N-d}$, respectively. We assume that $p_{N-d}(L)$ is dense, and $p_d$ is invertible on $p_d(L)$. Then the quasicrystal (called a {\it model set}) $\Gamma$ associated with $L$ and $W$ is defined to be $p_d\big(L\cap (\mathbb{R}^d\times W)\big)$. 

We put $\Lambda=p_d(L^*)$. It should be remarked that for each $\boldsymbol{\xi}\in \Lambda$, there exists a unique $\boldsymbol{\xi}'\in \mathbb{R}^{N-d}$ such that $(\boldsymbol{\xi},\boldsymbol{\xi}')\in L^*$. Indeed, if  $(\boldsymbol{\xi},\boldsymbol{\xi}'')\in L^*$, then $({\bf 0},\boldsymbol{\xi}'-\boldsymbol{\xi}'')\in L^*$, and hence $\mathbb{Z}\ni\langle ({\bf 0},\boldsymbol{\xi}'-\boldsymbol{\xi}''), \boldsymbol{\alpha}\rangle=\langle \boldsymbol{\xi}'-\boldsymbol{\xi}'', p_{N-d}(\boldsymbol{\alpha})\rangle$ for every $\boldsymbol{\alpha}\in L$. Since $p_{N-d}(L)$ is dense, we conclude that $\boldsymbol{\xi}'-\boldsymbol{\xi}''={\bf 0}$.

Let us write down a generalized Poisson formula for $\Gamma$ in a formal way.
Let $f$ be a compactly supported smooth function on $\mathbb{R}^d$, and let $\chi_W$ be the indicator function of the window $W\subset \mathbb{R}^{N-d}$. Define the compactly supported function $F$ on $\mathbb{R}^N$ by setting $F({\bf x},{\bf x}')=f({\bf x})\chi_W({\bf x}')$~ (${\bf x}\in \mathbb{R}^d, {\bf x}'\in \mathbb{R}^{N-d}$). Applying the Poisson summation formula to $F$, we obtain
$$
\sum_{{\bf z}\in \Gamma}f({\bf z})
=\sum_{\boldsymbol{\alpha}\in L}F(\boldsymbol{\alpha})={\rm vol}(D_L)^{-1}\sum_{\boldsymbol{\beta}\in L^*}\widehat{F}(\boldsymbol{\beta}),
$$  
which is, of course, a ``formal" identity because the right-hand side does not necessarily converge. Pretending that this is a genuine identity and noting 
$$
\widehat{F}(\boldsymbol{\beta})=\hat{f}(\boldsymbol{\xi})\widehat{\chi_W}(\boldsymbol{\xi}')\quad (\boldsymbol{\beta}=(\boldsymbol{\xi},\boldsymbol{\xi}')\in \mathbb{R}^d\times \mathbb{R}^{N-d}), 
$$
we get
\begin{equation}\label{eq:formal}
\sum_{{\bf z}\in \Gamma}f({\bf z})=\sum_{\boldsymbol{\xi}\in \Lambda}a(\boldsymbol{\xi})\hat{f}(\boldsymbol{\xi}), 
\end{equation}
where, for $\boldsymbol{\xi}\in \Lambda$, we put
$$
a(\boldsymbol{\xi})={\rm vol}(D_L)^{-1}\widehat{\chi_W}(\boldsymbol{\xi}') \quad((\boldsymbol{\xi},\boldsymbol{\xi}')\in L^*).
$$ 

We may justify Eq.~\!\ref{eq:formal} as follows. Let $U_{1/N}(W)$ be the $1/N$-neighborhood of $W$, and take a smooth function $g_N$ on $\mathbb{R}^{N-d}$ satisfying $0\leq g_N({\bf x}')\leq 1$ and
$$
g_N({\bf x}')=\begin{cases}
1 & ({\bf x}'\in W)\\
0 & ({\bf x}'\in U_{1/N}(W)^c)
\end{cases}.
$$
Put $F_N({\bf x},{\bf x}')=f({\bf x})g_N({\bf x}')$. 
If we take $N\gg 1$, we have $({\rm supp}~\!f\times U_{1/N}(W))\cap L=({\rm supp}~\!f\times W)\cap L$, so that, if $f({\bf z})g_N({\bf z}')\neq 0$ for $\boldsymbol{\alpha}=({\bf z},{\bf z}')\in L$, then $({\bf z},{\bf z}')\in ({\rm supp}~\!f\times U_{1/N}(W))\cap L=({\rm supp}~\!f\times W)\cap L$, and hence ${\bf z}\in \Gamma$ and $F_N(\boldsymbol{\alpha})=f({\bf z})$. We thus have
\begin{eqnarray*}
\sum_{{\bf z}\in \Gamma}f({\bf z})
&=&\sum_{\boldsymbol{\alpha}\in L}F_N(\boldsymbol{\alpha})={\rm vol}(D_L)^{-1}\sum_{\boldsymbol{\beta}\in L^*}\widehat{F_N}(\boldsymbol{\beta})\\
&=&{\rm vol}(D_L)^{-1}\sum_{(\boldsymbol{\xi},\boldsymbol{\xi}')\in L^*}\hat{f}(\boldsymbol{\xi})\widehat{g_N}(\boldsymbol{\xi}')\\
&=&\sum_{\boldsymbol{\xi}\in \Lambda}a_N(\boldsymbol{\xi})\hat{f}(\boldsymbol{\xi}),
\end{eqnarray*}
where $a_N(\boldsymbol{\xi})={\rm vol}(D_L)^{-1}\widehat{g_N}(\boldsymbol{\xi}')$. Obviously $\displaystyle\lim_{N\to\infty}a_N(\boldsymbol{\xi})=a(\boldsymbol{\xi})$.

A typical example of model sets is the set of nodes in a {\it Penrose tiling}\index{Penrose tiling} discovered by R. Penrose in 1973/1974, which is a remarkable non-periodic tiling generated by an aperiodic set of prototiles (see de Bruijn \cite{debr} for the proof of the fact that a Penrose tiling is obtained by the cut and projection method).

\vspace{-0.3cm}
\begin{figure}[htbp]
\begin{center}
\includegraphics[width=.4\linewidth]{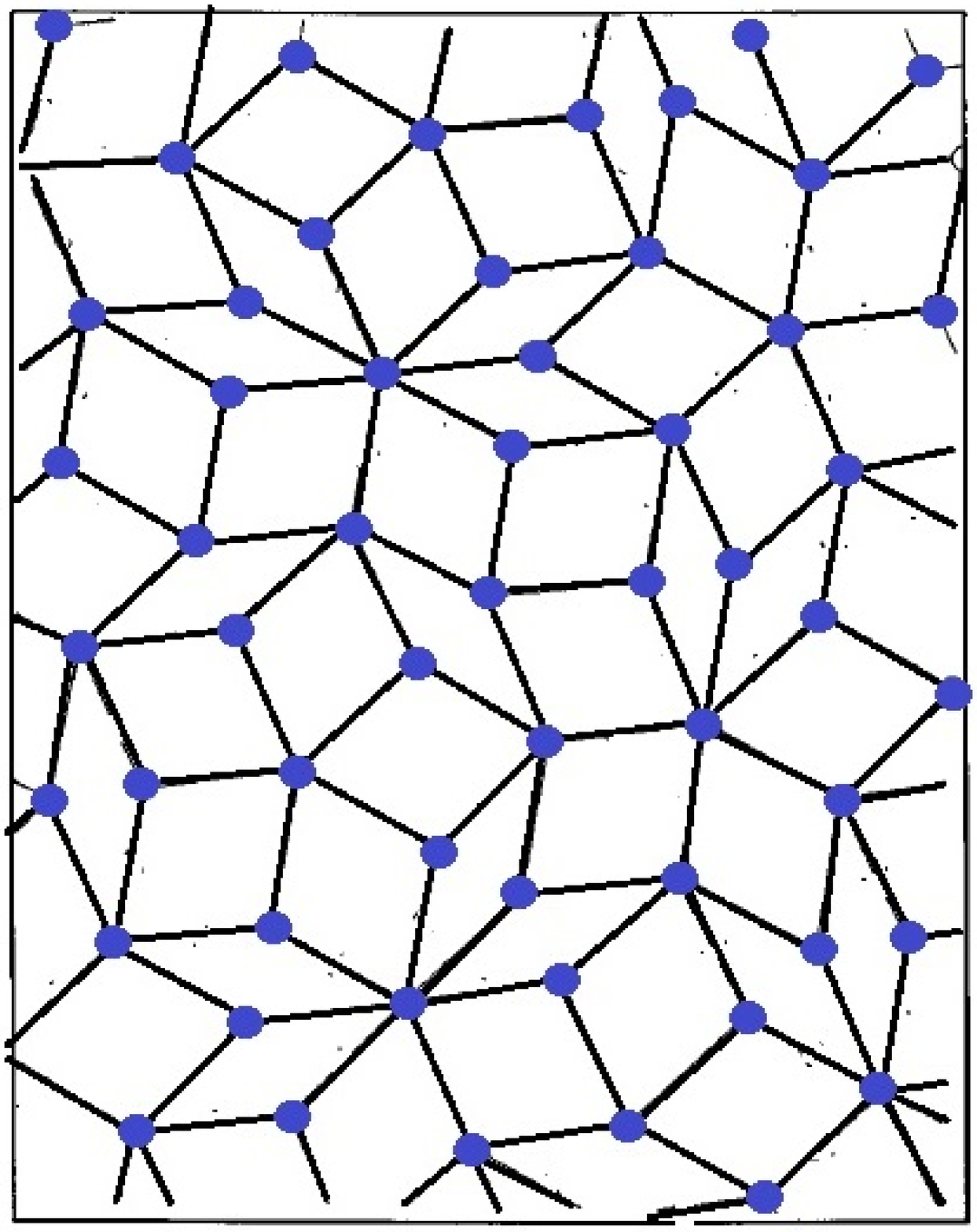}
\end{center}
\vspace{-1cm}
\caption{A Penrose tiling}\label{fig:penrose}
\end{figure}

\section{Is $\mathbb{Z}^d_{\rm prim}$ a quasicrystal?}
It is natural to ask whether $\mathbb{Z}^d_{\rm prim}$ is a quasicrystal. The answer is ``No." However $\mathbb{Z}^d_{\rm prim}$ is {\it nearly} a quasicrystal of Poisson type\index{near quasicrystal of Poisson type}.  

To see this, take a look again at the identity
\begin{equation*}\label{eq:vinoxx}
\sum_{{\bf z}\in \mathbb{Z}^d_{\rm prim}}f({\bf z})=\sum_{k=1}^{\infty}\mu(k)\sum_{{\bf w}\in \mathbb{Z}^d\backslash\{{\bf 0}\}}f(k{\bf w}).
\end{equation*}
Suppose that ${\rm supp}~\!f\subset B_N({\bf 0})$. Then applying the Poisson summation formula, we obtain
\begin{eqnarray*}
\sum_{{\bf z}\in \mathbb{Z}^d_{\rm prim}}f({\bf z})&=&\sum_{k=1}^{\infty}\mu(k)\sum_{{\bf w}\in \mathbb{Z}^d\backslash \{0\}}f(k{\bf w})\\
&=&\sum_{k=1}^N\mu(k)\Big[\sum_{{\bf w}\in \mathbb{Z}^d}f(k{\bf w})-f({\bf 0})\Big]\\
&=& \sum_{k=1}^N\mu(k)k^{-d}\sum_{\boldsymbol{\xi}\in k^{-1}\mathbb{Z}^d}\hat{f}(\boldsymbol{\xi})-\Big(\sum_{k=1}^N\mu(k)\Big)f({\bf 0}).
\end{eqnarray*}
Now for $\boldsymbol{\xi}\in \mathbb{Q}^d$, we write
$$
\boldsymbol{\xi}=\Big(\frac{b_1}{a_1},\ldots,\frac{b_d}{a_d}\Big), \quad {\rm gcd}(a_i,b_i)=1, a_i>0,
$$
and put
$
n(\boldsymbol{\xi})={\rm lcm}(a_1,\ldots,a_d)
$.
Then 
$
\boldsymbol{\xi}\in k^{-1}\mathbb{Z}^d~\Longleftrightarrow~n(\boldsymbol{\xi})|k
$, and hence
$$
\sum_{{\bf z}\in \mathbb{Z}^d_{\rm prim}}f({\bf z})=\sum_{k=1}^N\mu(k)k^{-d}\sum_{\boldsymbol{\xi}\in \mathbb{Q}^d\atop n(\boldsymbol{\xi})\mid k}\hat{f}(\boldsymbol{\xi})-\Big(\sum_{k=1}^N\mu(k)\Big)f({\bf 0}),
$$
where we should note that the first term in the right-hand side is an absolutely convergent series. 
To rewrite the right-hand side further, consider 
\begin{eqnarray*}
&&\mathbb{Q}^d_N=\{\boldsymbol{\xi}\in \mathbb{Q}^d|~\!n(\boldsymbol{\xi})\leq N\}, \\
&&A=\{(k,\boldsymbol{\xi})|~\!k=1,\ldots, N, ~\boldsymbol{\xi}\in \mathbb{Q}^d,~n(\boldsymbol{\xi})|k\},\\
&&B=\{(\ell,\boldsymbol{\xi})|~\!1\leq \ell \leq Nn(\boldsymbol{\xi}){}^{-1},~\boldsymbol{\xi}\in \mathbb{Q}^d_N\}.
\end{eqnarray*}
Then the map $(k,\boldsymbol{\xi})\mapsto (kn(\boldsymbol{\xi}){}^{-1},\boldsymbol{\xi})$ is a bijection of $A$ onto $B$. Therefore we get
$$
\sum_{{\bf z}\in \mathbb{Z}^d_{\rm prim}}f({\bf z})=\sum_{\boldsymbol{\xi}\in \mathbb{Q}^d_N}\sum_{1\leq \ell\leq N/n(\boldsymbol{\xi})}\frac{\mu(\ell n(\boldsymbol{\xi}))}{(\ell n(\boldsymbol{\xi}))^{d}}\hat{f}(\boldsymbol{\xi})-\Big(\sum_{k=1}^N\mu(k)\Big)f({\bf 0}).
$$
Clearly
$$
\mu(\ell n(\boldsymbol{\xi}))=\begin{cases}
\mu(\ell)\mu(n(\boldsymbol{\xi})) & ({\rm gcd}(\ell,n(\boldsymbol{\xi}))=1)\\
0 &({\rm gcd}(\ell,n(\boldsymbol{\xi}))>1).
\end{cases}
$$
Therefore putting
\begin{eqnarray*}
&&a_N(\boldsymbol{\xi})=\frac{\mu(n(\boldsymbol{\xi}))}{n(\boldsymbol{\xi}){}^{d}}\sum_{1\leq \ell\leq N/n(\boldsymbol{\xi})\atop {\rm gcd}(\ell,n(\boldsymbol{\xi}))=1}\frac{\mu(\ell)}{\ell^d},\\
&&\Lambda_N=\{\boldsymbol{\xi}\in \mathbb{Q}^d_N|~\!\mu(n(\boldsymbol{\xi}))\neq 0\},
\end{eqnarray*}
we get
\begin{equation*}\label{eq:lambda}
\sum_{{\bf z}\in \mathbb{Z}^d_{\rm prim}}f({\bf z})= \sum_{\boldsymbol{\xi}\in \Lambda_N}a_N(\boldsymbol{\xi})\hat{f}(\boldsymbol{\xi})-\Big(\sum_{k=1}^N\mu(k)\Big)f({\bf 0}).
\end{equation*}
Furthermore, if we put
\begin{eqnarray*}
&&\Lambda=\{\boldsymbol{\xi}\in \mathbb{Q}^d|~\!\mu(n(\boldsymbol{\xi}))\neq 0\},\\ 
&& a(\boldsymbol{\xi})=\frac{\mu(n_{\boldsymbol{\boldsymbol{\xi}}})}{n(\boldsymbol{\xi}){}^d}\zeta(d)^{-1}\prod_{p|n(\boldsymbol{\xi})}\big(1-p^{-d}\big)^{-1}\quad (\boldsymbol{\xi}\in \Lambda),
\end{eqnarray*}
then
\begin{eqnarray*}
\Lambda=\bigcup_{N=1}^{\infty}\Lambda_N,\qquad
\lim_{N\to\infty}a_N(\boldsymbol{\xi})=\frac{\mu(n(\boldsymbol{\xi}))}{n(\boldsymbol{\xi}){}^{d}}\sum_{\ell=1 \atop {\rm gcd}(\ell,n(\boldsymbol{\xi}))=1}^{\infty}\frac{\mu(\ell)}{\ell^d}=a(\boldsymbol{\xi}).
\end{eqnarray*}
This implies that if the ``extra term" $\displaystyle\Big(\sum_{k=1}^N\mu(k)\Big)f({\bf 0})$ is ignored, then the set $\mathbb{Z}^d_{\rm prim}$ looks like a quasicrystal of Poisson type. This is the reason why we say that $\mathbb{Z}^d_{\rm prim}$ is {\it nearly} a quasicrystal of Poisson type.

\begin{rem}{\rm 
(1) Applying Eq.~\!\ref{thm:pptx1}, we obtain
$$
\sum_{{\bf z}\in (\mathbb{Z}^{\rm odd})^2_{\rm prim}}f({\bf z})=\sum_{\boldsymbol{\xi}\in \mathbb{Q}_{2N}^2}
\Big(\sum_{k\geq 1, h\geq 0 \atop n(\boldsymbol{\xi})|k2^{h+1}}^{k2^h\leq N}\frac{\mu(k)}{k^2}\frac{1}{2^{2h+2}}e^{\pi i k2^{h+1}\langle\boldsymbol{\xi},{\bf 1}\rangle}
\Big)\hat{f}(\boldsymbol{\xi}),
$$
where ${\rm supp}~\!f\subset B_N({\bf 0})$ and ${\bf 1}=(1,1)$. This implies that $(\mathbb{Z}^{\rm odd})^2_{\rm prim}$ is a quasicrystal of Poisson type. The reason why no extra terms appear in this case is that  $(\mathbb{Z}^{\rm odd})^2=2\mathbb{Z}^2+{\bf 1}$ is a full lattice.

\medskip

(2) In much the same manner as above, we get
$$
\sum_{{\bf z}\in \mathbb{Z}^d_{\rm prim}}f({\bf z})e^{2\pi i \langle {\bf z},\boldsymbol{\eta}\rangle}
=\sum_{\boldsymbol{\xi}\in \Lambda_N}a_N(\boldsymbol{\xi})\hat{f}(\boldsymbol{\xi}-\boldsymbol{\eta})
-\Big(\sum_{k=1}^N\mu(k)\Big)f({\bf 0}).
$$
Using this identity, we can show
$$
\lim_{\epsilon \to +0}
\sigma(f^{\epsilon},\mathbb{Z}^d_{\rm prim},\omega_{\boldsymbol{\eta}})=a(\boldsymbol{\eta})\int_{\mathbb{R}^d}f({\bf x})~\!d{\bf x},
$$
that is, $(\mathbb{Z}^d_{\rm prim},\omega_{\boldsymbol{\eta}})$ has constant density for $\boldsymbol{\eta}$ with $\mu(n(\boldsymbol{\eta}))\neq 0$.

\medskip

(3) An interesting problem related to quasicrystals comes up in the study of
{\it non-trivial zeros}\index{non-trivial zeros of the Riemann zeta function} of the Riemann zeta function (thus we come across
another Riemann's work,\footnote{{\it \"{U}ber die Anzahl der Primzahlen unter einer gegebenen Gr\"{o}sse, 1859.}} which were to change the direction of mathematical
research in a most significant way).

We put
$$
\Gamma^{\rm zero}=\{{\rm Im}~\!s\in \mathbb{R}|~\zeta(s)=0,~0<{\rm Re}~\!s<1\}.
$$
Under the Riemann Hypothesis (RH), one may say that $\Gamma^{\rm zero}$ is {\it nearly} a quasicrystal of Poisson type of $1$-dimension (cf.~\!Dyson \!\cite{dy}). Actually a version of Riemann's explicit formula\index{Riemann's explicit formula} looks like a generalized Poisson formula (see Iwaniec and Kowalski \cite{iwa}):
\begin{eqnarray*}
&&\sum_{\rho}f\left(\frac{\rho-1/2}{i}\right)=
f\Big(\frac{1}{2i}\Big)+f\Big(-\frac{1}{2i}\Big)\nonumber\\
&&+\frac{1}{2\pi}\int_{-\infty}^{\infty}f(u){\rm Re}\frac{\Gamma'}{\Gamma}
\Big(\frac{1}{4}+\frac{iu}{2} \Big)du\nonumber\\
&&-\frac{1}{2\pi}\hat{f}(0)\log \pi-\frac{1}{2\pi}\sum_{m=1}^{\infty}\sum_{p}\frac{\log p}{p^{m/2}}
\left(\hat{f}\Big(\frac{\log p^m}{2\pi}\Big)+\hat{f}\Big(-\frac{\log p^m}{2\pi}\Big)
\right).\label{eq:explicit}
\end{eqnarray*}
where 
$\{\rho\}$ is the set of zeros of $\zeta(s)$ with $0<{\rm Re}~\!\rho<1$, 
$\displaystyle\sum_{p}$ is the sum over all primes, and $\Gamma'/\Gamma$ is the logarithmic derivative of the gamma function. Notice that, under the RH, the sum in the left-hand side is written as $\displaystyle\sum_{{\bf z}\in \Gamma^{\rm zero}}f({\bf z})$.\footnote{The {\it simple zero conjecture} says that all zeros $\rho$ are simple. In the case that we do not assume this conjecture, we think of $\Gamma^{\rm zero}$ as a weighted set with the weight $\omega(\rho)={\rm ord}_{\rho}(\zeta)$.} What we should stress here is that the test function $f(s)$ is not arbitrary, and is supposed to be analytic in the strip $|{\rm Im}~\!s| \leq1/2 + \epsilon$ for some $\epsilon > 0$, and to satisfy  $|f(s)| \leq (1 + |s|)^{-(1+\delta)}$ for some $\delta>0$ when $|{\rm Re}~\!s|\to\infty$. This restriction on $f$ together with the extra terms in the formula above says that $\Gamma^{\rm zero}$ is not a {\it genuine} quasicrystal of Poisson type. Furthermore $\Gamma^{\rm zero}$ does not have the Delone property. \hfill$\Box$

}
\end{rem}


\end{document}